\documentclass{amsart}
\usepackage{amssymb}
\usepackage{amsmath}
\usepackage{amsfonts}

\newtheorem{theorem}{Theorem}
\theoremstyle{plain}

\newtheorem{corollary}{Corollary}

\newtheorem{definition}{Definition}

\newtheorem{proposition}{Proposition}

\numberwithin{equation}{section}
\input{tcilatex}

\begin{document}
\title[weighted $q$-Hardy-littlewood-type maximal operator ]{A note on the
weighted $q$-Hardy-littlewood-type maximal operator with respect to $q$%
-Volkenborn integral in the $p$-adic integer ring}
\author{Serkan Araci}
\address{University of Gaziantep, Faculty of Science and Arts, Department of
Mathematics, 27310 Gaziantep, TURKEY}
\email{mtsrkn@hotmail.com}
\author{Mehmet Acikgoz}
\address{University of Gaziantep, Faculty of Science and Arts, Department of
Mathematics, 27310 Gaziantep, TURKEY}
\email{acikgoz@gantep.edu.tr}
\subjclass[2000]{Primary 05A10, 11B65; Secondary 11B68, 11B73.}
\keywords{$q$-Volkenborn integral, Hardy-littlewood theorem, $p$-adic
analysis, $q$-analysis}

\begin{abstract}
The essential aim of this paper is to define weighted $q$%
-Hardy-littlewood-type maximal operator by means of $p$-adic $q$-invariant
distribution on $%
\mathbb{Z}
_{p}$. Moreover, we give some interesting properties concerning this type
maximal operator.
\end{abstract}

\maketitle

\section{Introduction, Definitions and Notations}

\bigskip Recently, $q$-analysis has served as a structure between
mathematics and physics. Therefore, there is a significant increase of
activity in the area of the $q$-analysis due to applications of the $q$%
-analysis in mathematics, statistics and physics.

$p$-adic numbers also play a vital and important role in mathematics. $p$%
-adic numbers were invented by the German mathematician Kurt Hensel \cite%
{Hensel}, around the end of the nineteenth century. In spite of their being
already one hundred years old, these numbers are still today enveloped in an
aura of mystery within the scientific community.

The $p$-adic $q$-integral (or $q$-Volkenborn integral) are originally
constructed by Kim \cite{Kim 4}. The $q$-Volkenborn integral is used in
Mathematical Physics for example the functional equation of the $q$-Zeta
function, the $q$-Stirling numbers and $q$-Mahler theory of integration with
respect to the ring $%
\mathbb{Z}
_{p}$ together with Iwasawa's $p$-adic $q$-$L$ function.

T. Kim, by using $q$-Volkenborn integral, introduced a novel
Lebesgue-Radon-Nikodym type theorem. He is given some interesting properties
concerning Lebesgue-Radon-Nikodym theorem. After, Jang \cite{Jang} defined $%
q $-extension of Hardy-Littlewood-type maximal operator by means of $q$%
-Volkenborn integral on $%
\mathbb{Z}
_{p}.$ Next, Kim, Choi and Kim \cite{Kim 7} defined weighted
Lebesgue-Radon-Nikodym theorem. They also gave some interesting properties
of this type theorem.

By the same motivation of the above studies, in this paper, we construct
weighted $q$-Hardy-littlewood-type maximal operator by the means of $p$-adic 
$q$-integral on $%
\mathbb{Z}
_{p}$. We also give some interesting properties of this type operator.

Imagine that $p$ be a fixed prime number. Let $%
\mathbb{Q}
_{p}$ be the field of $p$-adic rational numbers and let $%
\mathbb{C}
_{p}$ be the completion of algebraic closure of $%
\mathbb{Q}
_{p}$.

Thus, 
\begin{equation*}
\mathbb{Q}
_{p}=\left\{ x=\sum_{n=-k}^{\infty }a_{n}p^{n}:0\leq a_{n}\leq p-1\right\} .
\end{equation*}

Then $%
\mathbb{Z}
_{p}$ is an integral domain, which is defined by 
\begin{equation*}
\mathbb{Z}
_{p}=\left\{ x=\sum_{n=0}^{\infty }a_{n}p^{n}:0\leq a_{n}\leq p-1\right\} ,
\end{equation*}

or 
\begin{equation*}
\mathbb{Z}
_{p}=\left\{ x\in 
\mathbb{Q}
_{p}:\left\vert x\right\vert _{p}\leq 1\right\} .
\end{equation*}

In this paper, we assume that $q\in 
\mathbb{C}
_{p}$ with $\left\vert 1-q\right\vert _{p}<1$ as an indeterminate.

The $p$-adic absolute value $\left\vert .\right\vert _{p}$, is normally
defined by 
\begin{equation*}
\left\vert x\right\vert _{p}=\frac{1}{p^{r}}\text{,}
\end{equation*}

where $x=p^{r}\frac{s}{t}$ with $\left( p,s\right) =\left( p,t\right)
=\left( s,t\right) =1$ and $r\in 
\mathbb{Q}
$.

A $p$-adic Banach space $B$ is a $%
\mathbb{Q}
_{p}$-vector space with a lattice $B^{0}$ ($%
\mathbb{Z}
_{p}$-module) separated and complete for $p$-adic topology, ie., 
\begin{equation*}
B^{0}\simeq \lim_{\overleftarrow{n\in 
\mathbb{N}
}}B^{0}/p^{n}B^{0}\text{.}
\end{equation*}

For all $x\in B$, there exists $n\in 
\mathbb{Z}
$, such that $x\in p^{n}B^{0}$. Define 
\begin{equation*}
v_{B}\left( x\right) =\sup_{n\in 
\mathbb{N}
\cup \left\{ +\infty \right\} }\left\{ n:x\in p^{n}B^{0}\right\} \text{.}
\end{equation*}

It satisfies the following properties:%
\begin{eqnarray*}
v_{B}\left( x+y\right) &\geq &\min \left( v_{B}\left( x\right) ,v_{B}\left(
y\right) \right) \text{,} \\
v_{B}\left( \beta x\right) &=&v_{p}\left( \beta \right) +v_{B}\left(
x\right) \text{, if }\beta \in 
\mathbb{Q}
_{p}\text{.}
\end{eqnarray*}

Then, $\left\Vert x\right\Vert _{B}=p^{-v_{B}\left( x\right) }$ defines a
norm on $B,$ such that $B$ is complete for $\left\Vert .\right\Vert _{B}$
and $B^{0}$ is the unit ball.

A measure on $%
\mathbb{Z}
_{p}$ with values in a $p$-adic Banach space $B$ is a continuous linear map%
\begin{equation*}
f\mapsto \int f\left( x\right) \mu =\int_{%
\mathbb{Z}
_{p}}f\left( x\right) \mu \left( x\right)
\end{equation*}

from $C^{0}\left( 
\mathbb{Z}
_{p},%
\mathbb{C}
_{p}\right) $, (continuous function on $%
\mathbb{Z}
_{p}$) to $B$. We know that the set of locally constant functions from $%
\mathbb{Z}
_{p}$ to $%
\mathbb{Q}
_{p}$ is dense in $C^{0}\left( 
\mathbb{Z}
_{p},%
\mathbb{C}
_{p}\right) $ so.

Explicitly, for all $f\in C^{0}\left( 
\mathbb{Z}
_{p},%
\mathbb{C}
_{p}\right) $, the locally constant functions 
\begin{equation*}
f_{n}=\sum_{i=0}^{p^{n}-1}f\left( i\right) 1_{i+p^{n}%
\mathbb{Z}
_{p}}\rightarrow \text{ }f\text{ in }C^{0}\text{.}
\end{equation*}

Now if ~$\mu \in \boldsymbol{D}_{0}\left( 
\mathbb{Z}
_{p},%
\mathbb{Q}
_{p}\right) $, set $\mu \left( i+p^{n}%
\mathbb{Z}
_{p}\right) =\int_{%
\mathbb{Z}
_{p}}1_{i+p^{n}%
\mathbb{Z}
_{p}}\mu $. Then $\int_{%
\mathbb{Z}
_{p}}f\mu $ is given by the following \textquotedblleft Riemann
sums\textquotedblright 
\begin{equation*}
\int_{%
\mathbb{Z}
_{p}}f\mu =\lim_{n\rightarrow \infty }\sum_{i=0}^{p^{n}-1}f\left( i\right)
\mu \left( i+p^{n}%
\mathbb{Z}
_{p}\right) \text{.}
\end{equation*}

T. Kim defined $\mu _{q}$ as follows:%
\begin{equation*}
\mu _{q}\left( \xi +dp^{n}%
\mathbb{Z}
_{p}\right) =\frac{q^{\xi }}{\left[ dp^{n}\right] _{q}}
\end{equation*}

and this can be extended to a distribution on $%
\mathbb{Z}
_{p}$. This distribution yields an integral in the case $d=1$.

So, $q$-Volkenborn integral was defined by T. Kim as follows:%
\begin{eqnarray}
&&I_{q}\left( f\right) =\int_{%
\mathbb{Z}
_{p}}f\left( \xi \right) d\mu _{q}\left( \xi \right)  \label{equation 5} \\
&=&\lim_{n\rightarrow \infty }\frac{1}{\left[ p^{n}\right] _{q}}\sum_{\xi
=0}^{p^{n}-1}f\left( \xi \right) q^{\xi }\text{, (for details, see \cite{Kim
4}, \cite{Kim 5}). }  \notag
\end{eqnarray}

Where $\left[ x\right] _{q}$ is a $q$-extension of $x$ which is defined by%
\begin{equation*}
\left[ x\right] _{q}=\frac{1-q^{x}}{1-q}\text{,}
\end{equation*}

note that $\lim_{q\rightarrow 1}\left[ x\right] _{q}=x$ cf. \cite{Kim 2}, 
\cite{Kim 3}, \cite{Kim 4}, \cite{Kim 5}, \cite{Jang}.

Let $d$ be a fixed positive integer with $\left( p,d\right) =1$. We now set%
\begin{eqnarray*}
X &=&X_{d}=\lim_{\overleftarrow{n}}%
\mathbb{Z}
/dp^{n}%
\mathbb{Z}
, \\
X_{1} &=&%
\mathbb{Z}
_{p}, \\
X^{\ast } &=&\underset{\underset{\left( a,p\right) =1}{0<a<dp}}{\cup }a+dp%
\mathbb{Z}
_{p}, \\
a+dp^{n}%
\mathbb{Z}
_{p} &=&\left\{ x\in X\mid x\equiv a\left( \func{mod}p^{n}\right) \right\} ,
\end{eqnarray*}

where $a\in 
\mathbb{Z}
$ satisfies the condition $0\leq a<dp^{n}$. For $f\in UD\left( 
\mathbb{Z}
_{p},%
\mathbb{C}
_{p}\right) $,%
\begin{equation*}
\int_{%
\mathbb{Z}
_{p}}f\left( x\right) d\mu _{q}\left( x\right) =\int_{X}f\left( x\right)
d\mu _{q}\left( x\right) ,
\end{equation*}

(for details, see \cite{Kim 8}).

By the meaning of $q$-Volkenborn integral, we consider below strongly $p$%
-adic $q$-invariant distribution $\mu _{q}$ on $%
\mathbb{Z}
_{p}$ in the form 
\begin{equation*}
\left\vert \left[ p^{n}\right] _{q}\mu _{q}\left( a+p^{n}%
\mathbb{Z}
_{p}\right) -\left[ p^{n+1}\right] _{q}\mu _{q}\left( a+p^{n+1}%
\mathbb{Z}
_{p}\right) \right\vert <\delta _{n},
\end{equation*}

where $\delta _{n}\rightarrow 0$ as $n\rightarrow \infty $ and $\delta _{n}$
is independent of $a$. Let $f\in UD\left( 
\mathbb{Z}
_{p},%
\mathbb{C}
_{p}\right) $, for any $a\in 
\mathbb{Z}
_{p}$, we assume that the weight function $\omega \left( x\right) $ is
defined by $\omega \left( x\right) =\omega ^{x}$ where $\omega \in 
\mathbb{C}
_{p}$ with $\left\vert 1-\omega \right\vert _{p}<1$. We define the weighted
measure on $%
\mathbb{Z}
_{p}$ as follows:%
\begin{equation}
\mu _{f,q}^{\left( \omega \right) }\left( a+p^{n}%
\mathbb{Z}
_{p}\right) =\int_{a+p^{n}%
\mathbb{Z}
_{p}}\omega ^{\xi }f\left( \xi \right) d\mu _{q}\left( \xi \right)
\label{equation 2}
\end{equation}

where the integral is the $q$-Volkenborn integral. By (\ref{equation 2}), we
easily note that $\mu _{f,q}^{\left( \omega \right) }$ is a strongly
weighted measure on $%
\mathbb{Z}
_{p}$. That is,%
\begin{eqnarray*}
&&\left\vert \left[ p^{n}\right] _{q}\mu _{f,q}^{\left( \omega \right)
}\left( a+p^{n}%
\mathbb{Z}
_{p}\right) -\left[ p^{n+1}\right] _{q}\mu _{f,q}^{\left( \omega \right)
}\left( a+p^{n+1}%
\mathbb{Z}
_{p}\right) \right\vert _{p} \\
&=&\left\vert \sum_{x=0}^{p^{n}-1}\omega ^{x}f\left( x\right)
q^{x}-\sum_{x=0}^{p^{n}}\omega ^{x}f\left( x\right) q^{x}\right\vert _{p} \\
&\leq &\left\vert \frac{f\left( p^{n}\right) \omega ^{p^{n}}q^{p^{n}}}{p^{n}}%
\right\vert _{p}\left\vert p^{n}\right\vert _{p} \\
&\leq &Cp^{-n}
\end{eqnarray*}

Thus, we get the following proposition.

\begin{proposition}
For $f,g\in UD\left( 
\mathbb{Z}
_{p},%
\mathbb{C}
_{p}\right) $, we have 
\begin{equation*}
\mu _{\alpha f+\beta g,q}^{\left( \omega \right) }\left( a+p^{n}%
\mathbb{Z}
_{p}\right) =\alpha \mu _{f,q}^{\left( \omega \right) }\left( a+p^{n}%
\mathbb{Z}
_{p}\right) +\beta \mu _{g,q}^{\left( \omega \right) }\left( a+p^{n}%
\mathbb{Z}
_{p}\right) \text{.}
\end{equation*}%
where $\alpha ,\beta $ are positive constants. Moreover,%
\begin{equation*}
\left\vert \left[ p^{n}\right] _{q}\mu _{f,q}^{\left( \omega \right) }\left(
a+p^{n}%
\mathbb{Z}
_{p}\right) -\left[ p^{n+1}\right] _{q}\mu _{f,q}^{\left( \omega \right)
}\left( a+p^{n+1}%
\mathbb{Z}
_{p}\right) \right\vert \leq Cp^{-n}
\end{equation*}%
where $C$ is positive constant.
\end{proposition}

Let $UD\left( 
\mathbb{Z}
_{p},%
\mathbb{C}
_{p}\right) $ be the space of uniformly differentiable functions on $%
\mathbb{Z}
_{p}$ with supnorm 
\begin{equation*}
\left\Vert f\right\Vert _{\infty }=\underset{x\in 
\mathbb{Z}
_{p}}{\sup }\left\vert f\left( x\right) \right\vert _{p}.
\end{equation*}

The difference quotient $\Delta _{1}f$ of $f$ is the function of two
variables given by 
\begin{equation*}
\Delta _{1}f\left( m,x\right) =\frac{f\left( x+m\right) -f\left( x\right) }{m%
},\text{ for all }x\text{, }m\in 
\mathbb{Z}
_{p}\text{, }m\neq 0\text{.}
\end{equation*}

A function $f:%
\mathbb{Z}
_{p}\rightarrow 
\mathbb{C}
_{p}$ is said to be a Lipschitz function if there exists a constant $M>0$ $%
\left( \text{the Lipschitz constant of }f\right) $ such that%
\begin{equation*}
\left\vert \Delta _{1}f\left( m,x\right) \right\vert \leq M\text{ for all }%
m\in 
\mathbb{Z}
_{p}\backslash \left\{ 0\right\} \text{ and }x\in 
\mathbb{Z}
_{p}.
\end{equation*}

The $%
\mathbb{C}
_{p}$ linear space consisting of all Lipschitz function is denoted by $%
Lip\left( 
\mathbb{Z}
_{p},%
\mathbb{C}
_{p}\right) $. This space is a Banach space with the respect to the norm $%
\left\Vert f\right\Vert _{1}=\left\Vert f\right\Vert _{\infty }\tbigvee
\left\Vert \Delta _{1}f\right\Vert _{\infty }$ (for more informations, see 
\cite{Kim 1}, \cite{Kim 2}, \cite{Kim 3}, \cite{Kim 4}, \cite{Kim 5}, \cite%
{Kim 6}, \cite{Jang}). The main aim of this paper is to define weighted $q$%
-extension of Hardy Littlewood type maximal operator. Moreover, we show the
boundedness of the weighted $q$-Hardy-littlewood-type maximal operator in
the $p$-adic integer ring.

\section{\qquad The weighted $q$-Hardy-littlewood-type maximal operator}

In view of (\ref{equation 2}) and the definition of $p$-adic $q$-integral on 
$%
\mathbb{Z}
_{p}$, we now start with the following theorem.

\begin{theorem}
Let $\mu _{q}^{\left( w\right) }$ be a strongly $p$-adic $q$-invariant in
the $p$-adic integer ring and $f\in UD\left( 
\mathbb{Z}
_{p},%
\mathbb{C}
_{p}\right) $. Then for any $n\in 
\mathbb{Z}
$ and any $\xi \in 
\mathbb{Z}
_{p}$, we have
\end{theorem}

$(1)$ $\int_{a+p^{n}%
\mathbb{Z}
_{p}}\omega ^{\frac{\xi }{p^{n}}}f\left( \xi \right) q^{-\frac{\xi }{p^{n}}%
}d\mu _{q^{p^{-n}}}\left( \xi \right) =\frac{\omega ^{\frac{a}{p^{n}}}}{%
\left[ p^{n}\right] _{q^{p^{-n}}}}\int_{%
\mathbb{Z}
_{p}}\omega ^{\xi }f\left( a+p^{n}\xi \right) q^{-\xi }d\mu _{q}\left( \xi
\right) ,$

$(2)$ $\int_{a+p^{n}%
\mathbb{Z}
_{p}}\omega ^{\frac{\xi }{p^{n}}}d\mu _{q^{p^{-n}}}\left( \xi \right) =\frac{%
\left( 1-q\right) \omega ^{\frac{a}{p^{n}}}q^{\frac{a}{p^{n}}}}{\left(
1-\omega q\right) \left[ p^{n}\right] _{q^{p^{-n}}}}\left( 1+\frac{\log
\omega }{\log q}\right) .$

\begin{proof}
(1) We see use of (\ref{equation 5}) and (\ref{equation 2})%
\begin{eqnarray*}
&&\int_{a+p^{n}%
\mathbb{Z}
_{p}}\omega ^{\frac{\xi }{p^{n}}}f\left( \xi \right) q^{-\frac{\xi }{p^{n}}%
}d\mu _{q^{p^{-n}}}\left( \xi \right) \\
&=&\lim_{m\rightarrow \infty }\frac{1}{\left[ p^{m+n}\right] _{q^{p^{-n}}}}%
\sum_{\xi =0}^{p^{m}-1}\omega ^{\frac{a+p^{n}\xi }{p^{n}}}f\left( a+p^{n}\xi
\right) q^{-\frac{a+p^{n}\xi }{p^{n}}}q^{\frac{a+p^{n}\xi }{p^{n}}} \\
&=&\omega ^{\frac{a}{p^{n}}}\lim_{m\rightarrow \infty }\frac{1}{\left[ p^{n}%
\right] _{q^{p^{-n}}}\left[ p^{m}\right] _{q}}\sum_{\xi =0}^{p^{m}-1}\omega
^{\xi }q^{-\xi }f\left( a+p^{n}\xi \right) q^{\xi } \\
&=&\frac{\omega ^{\frac{a}{p^{n}}}}{\left[ p^{n}\right] _{q^{p^{-n}}}}\int_{%
\mathbb{Z}
_{p}}\omega ^{\xi }f\left( a+p^{n}\xi \right) q^{-\xi }d\mu _{q}\left( \xi
\right) .
\end{eqnarray*}

(2) By the same method of (1), we easily see the following%
\begin{eqnarray*}
&&\int_{a+p^{n}%
\mathbb{Z}
_{p}}\omega ^{\frac{\xi }{p^{n}}}d\mu _{q^{p^{-n}}}\left( \xi \right) \\
&=&\lim_{m\rightarrow \infty }\frac{1}{\left[ p^{m+n}\right] _{q^{p^{-n}}}}%
\sum_{\xi =0}^{p^{m}-1}\omega ^{\frac{a+\xi p^{n}}{p^{n}}}q^{\frac{a+\xi
p^{n}}{p^{n}}} \\
&=&\frac{\omega ^{\frac{a}{p^{n}}}q^{\frac{a}{p^{n}}}}{\left[ p^{n}\right]
_{q^{p^{-n}}}}\lim_{m\rightarrow \infty }\frac{1}{\left[ p^{m}\right] _{q}}%
\sum_{\xi =0}^{p^{m}-1}\omega ^{\xi }q^{\xi } \\
&=&\frac{\omega ^{\frac{a}{p^{n}}}q^{\frac{a}{p^{n}}}\left( 1-q\right) }{%
\left[ p^{n}\right] _{q^{p^{-n}}}\left( 1-\omega q\right) }%
\lim_{m\rightarrow \infty }\frac{1-\left( \omega q\right) ^{p^{m}}}{%
1-q^{p^{m}}} \\
&=&\frac{\omega ^{\frac{a}{p^{n}}}q^{\frac{a}{p^{n}}}\left( 1-q\right) }{%
\left[ p^{n}\right] _{q^{p^{-n}}}\left( 1-\omega q\right) }\left( 1+\frac{%
\log \omega }{\log q}\right)
\end{eqnarray*}

Since $\underset{m\rightarrow \infty }{\lim }q^{p^{m}}=1$ for $\left\vert
1-q\right\vert _{p}<1,$ our assertion follows.
\end{proof}

Now, we are ready to give definition of weighted $q$-extension of
Hardy-littlewood-type maximal operator related to $p$-adic $q$-integral on $%
\mathbb{Z}
_{p}$ with a strong $p$-adic $q$-invariant distribution $\mu _{q}$ in the $p$%
-adic integer ring.

\begin{definition}
Let $\mu _{q}^{\left( \omega \right) }$ be a strongly $p$-adic $q$-invariant
distribution in the $p$-adic integer ring and $f\in UD\left( 
\mathbb{Z}
_{p},%
\mathbb{C}
_{p}\right) $. Then weighted $q$-extension of the Hardy-littlewood-type
maximal operator with respect to $p$-adic $q$-integral on $a+p^{n}%
\mathbb{Z}
_{p}$ is defined by 
\begin{equation*}
M_{p,q}^{\left( \omega \right) }f\left( a\right) =\underset{n\in 
\mathbb{Z}
}{\sup }\frac{1}{\mu _{1,q^{p^{-n}}}^{\left( w^{p^{-n}}\right) }\left(
x+p^{n}%
\mathbb{Z}
_{p}\right) }\int_{a+p^{n}%
\mathbb{Z}
_{p}}\omega ^{\frac{x}{p^{n}}}q^{-\frac{x}{p^{n}}}f\left( x\right) d\mu _{q^{%
\frac{1}{p^{n}}}}\left( x\right)
\end{equation*}%
for all $a\in 
\mathbb{Z}
_{p}$.
\end{definition}

We recall that famous Hardy-littlewood maximal operator $M_{\mu }$ is
defined by 
\begin{equation}
M_{\mu }f\left( a\right) =\underset{a\in Q}{\sup }\frac{1}{\mu \left(
Q\right) }\int_{Q}\left\vert f\left( x\right) \right\vert d\mu \left(
x\right) ,  \label{equation 3}
\end{equation}

where $f:%
\mathbb{R}
^{k}\rightarrow 
\mathbb{R}
^{k}$ is a locally bounded Lebesgue measurable function, $\mu $ is a
Lebesgue measure on $\left( -\infty ,\infty \right) $ and the supremum is
taken over all cubes $Q$ which are parallel to the coordinate axes. Note
that the boundedness of the Hardy-Littlewood maximal operator serves as one
of the most important tools used in the investigation of the properties of
variable exponent spaces (see \cite{Jang}). The essential aim of Theorem 1
is to deal with the weighted $q$-extension of the classical Hardy-Littlewood
maximal operator in the space of $p$-adic Lipschitz functions on $%
\mathbb{Z}
_{p}$ and to find the boundedness of them. By the meaning of Definition 1,
we get the following theorem.

\begin{theorem}
Let $f\in UD\left( 
\mathbb{Z}
_{p},%
\mathbb{C}
_{p}\right) $ and $x\in 
\mathbb{Z}
_{p}$, we get
\end{theorem}

\textit{(1)} $M_{p,q}^{\left( \omega \right) }f\left( a\right) =\frac{%
1-\omega q}{1-q}\frac{\log q}{\log \left( \omega q\right) }\underset{n\in 
\mathbb{Z}
}{\sup }q^{-\frac{x}{p^{n}}}\int_{%
\mathbb{Z}
_{p}}\omega ^{\xi }f\left( x+p^{n}\xi \right) q^{-\xi }d\mu _{q}\left( \xi
\right) $,

\textit{(2)} $\left\vert M_{p,q}^{\left( \omega \right) }f\left( a\right)
\right\vert _{p}\leq \underset{n\in 
\mathbb{Z}
}{\sup }\left\vert \frac{1-wq}{q^{\frac{x}{p^{n}}}}\frac{\log q}{\log \left(
\omega q\right) }\right\vert _{p}\left\Vert f\right\Vert _{1}\left\Vert
\left( \frac{q}{\omega }\right) ^{-\left( .\right) }\right\Vert _{L^{1}}$,

\textit{where} $\left\Vert \left( \frac{q}{\omega }\right) ^{-\left(
.\right) }\right\Vert _{L^{1}}=\int_{%
\mathbb{Z}
_{p}}\left( \frac{q}{\omega }\right) ^{-\xi }d\mu _{q}\left( \xi \right) $.

\begin{proof}
(1) Because of Theorem 1 and Definition 1, we see%
\begin{eqnarray*}
M_{p,q}^{\left( \omega \right) }f\left( a\right) &=&\underset{n\in 
\mathbb{Z}
}{\sup }\frac{1}{\mu _{1,q^{p^{-n}}}^{\left( w^{p^{-n}}\right) }\left(
x+p^{n}%
\mathbb{Z}
_{p}\right) }\int_{a+p^{n}%
\mathbb{Z}
_{p}}\omega ^{\frac{x}{p^{n}}}q^{-\frac{x}{p^{n}}}f\left( x\right) d\mu _{q^{%
\frac{1}{p^{n}}}}\left( x\right) \\
&=&\underset{n\in 
\mathbb{Z}
}{\sup }\frac{\left( 1-\omega q\right) \left[ p^{n}\right] _{q^{p^{-n}}}\log
q}{\left( 1-q\right) \omega ^{\frac{x}{p^{n}}}q^{\frac{x}{p^{n}}}\log \left(
\omega q\right) }\frac{\omega ^{\frac{x}{p^{n}}}}{\left[ p^{n}\right]
_{q^{p^{-n}}}}\int_{%
\mathbb{Z}
_{p}}\omega ^{\xi }f\left( a+p^{n}\xi \right) q^{-\xi }d\mu _{q}\left( \xi
\right) \\
&=&\frac{\left( 1-\omega q\right) \log q}{\left( 1-q\right) \log \left(
\omega q\right) }\underset{r\in 
\mathbb{Z}
}{\sup }\frac{1}{q^{\frac{x}{p^{n}}}}\int_{%
\mathbb{Z}
_{p}}\omega ^{\xi }f\left( a+p^{n}\xi \right) q^{-\xi }d\mu _{q}\left( \xi
\right)
\end{eqnarray*}

(2) On account of (1), we can derive the following%
\begin{eqnarray*}
\left\vert M_{p,q}^{\left( \omega \right) }f\left( a\right) \right\vert _{p}
&=&\left\vert \frac{\left( 1-\omega q\right) \log q}{\left( 1-q\right) \log
\left( \omega q\right) }\underset{n\in 
\mathbb{Z}
}{\sup }\frac{1}{q^{\frac{x}{p^{n}}}}\int_{%
\mathbb{Z}
_{p}}\omega ^{\xi }f\left( a+p^{n}\xi \right) q^{-\xi }d\mu _{q}\left( \xi
\right) \right\vert _{p} \\
&\leq &\left\vert \frac{\left( 1-\omega q\right) \log q}{\left( 1-q\right)
\log \left( \omega q\right) }\right\vert _{p}\underset{n\in 
\mathbb{Z}
}{\sup }\left\vert q^{-\frac{x}{p^{n}}}\int_{%
\mathbb{Z}
_{p}}f\left( a+p^{n}\xi \right) \left( \frac{q}{\omega }\right) ^{-\xi }d\mu
_{q}\left( \xi \right) \right\vert _{p} \\
&\leq &\left\vert \frac{\left( 1-\omega q\right) \log q}{\left( 1-q\right)
\log \left( \omega q\right) }\right\vert _{p}\underset{n\in 
\mathbb{Z}
}{\sup }\left\vert q^{-\frac{x}{p^{n}}}\right\vert _{p}\int_{%
\mathbb{Z}
_{p}}\left\vert f\left( a+p^{n}\xi \right) \right\vert _{p}\left\vert \left( 
\frac{q}{\omega }\right) ^{-\xi }\right\vert _{p}d\mu _{q}\left( \xi \right)
\\
&\leq &\left\vert \frac{\left( 1-\omega q\right) \log q}{\left( 1-q\right)
\log \left( \omega q\right) }\right\vert _{p}\underset{n\in 
\mathbb{Z}
}{\sup }\left\vert q^{-\frac{x}{p^{n}}}\right\vert _{p}\left\Vert
f\right\Vert _{1}\int_{%
\mathbb{Z}
_{p}}\left\vert \left( \frac{q}{\omega }\right) ^{-\xi }\right\vert _{p}d\mu
_{q}\left( \xi \right) \\
&=&\left\vert \frac{\left( 1-\omega q\right) \log q}{\left( 1-q\right) \log
\left( \omega q\right) }\right\vert _{p}\underset{n\in 
\mathbb{Z}
}{\sup }\left\vert q^{-\frac{x}{p^{n}}}\right\vert _{p}\left\Vert
f\right\Vert _{1}\left\Vert \left( \frac{q}{\omega }\right) ^{-\left(
.\right) }\right\Vert _{L^{1}}.
\end{eqnarray*}

Thus, we complete the proof of theorem.
\end{proof}

We note that Theorem 2 (2) shows the supnorm-inequality for the weighted $q$%
-Hardy-Littlewood-type maximal operator in the $p$-adic integer ring, in a
word, Theorem 2 (2) shows the following inequality%
\begin{equation}
\left\Vert M_{p,q}^{\left( \omega \right) }f\right\Vert _{\infty }=\underset{%
x\in 
\mathbb{Z}
_{p}}{\sup }\left\vert M_{p,q}^{\left( \omega \right) }f\left( x\right)
\right\vert _{p}\leq K\left\Vert f\right\Vert _{1}\left\Vert \left( \frac{q}{%
\omega }\right) ^{-\left( .\right) }\right\Vert _{L^{1}}  \label{equation 4}
\end{equation}

where $K=\left\vert \frac{\left( 1-\omega q\right) \log q}{\left( 1-q\right)
\log \left( \omega q\right) }\right\vert _{p}\underset{n\in 
\mathbb{Z}
}{\sup }\left\vert q^{-\frac{x}{p^{n}}}\right\vert _{p}$. By the equation (%
\ref{equation 4}), we get the following Corollary, which is the boundedness
for weighted $q$-Hardy-Littlewood-type maximal operator in the $p$-adic
integer ring.

\begin{corollary}
$M_{p,q}^{\left( \omega \right) }$ is a bounded operator from $UD\left( 
\mathbb{Z}
_{p},%
\mathbb{C}
_{p}\right) $ into $L^{\infty }\left( 
\mathbb{Z}
_{p},%
\mathbb{C}
_{p}\right) $, where $L^{\infty }\left( 
\mathbb{Z}
_{p},%
\mathbb{C}
_{p}\right) $ is the space of all $p$-adic supnorm-bounded functions with
the 
\begin{equation*}
\left\Vert f\right\Vert _{\infty }=\underset{x\in 
\mathbb{Z}
_{p}}{\sup }\left\vert f\left( x\right) \right\vert _{p}\text{,}
\end{equation*}%
for all $f\in L^{\infty }\left( 
\mathbb{Z}
_{p},%
\mathbb{C}
_{p}\right) $.
\end{corollary}

\end{document}